\documentstyle{amsart}
\title{Blow-up formulas for $(-2)$-spheres}
\author{Rogier Brussee}
\date{dg-ga/9412004}
\address{Fakult\"at f\"ur Mathematik Universit\"at Bielefeld \\
Postfach 100131\\
33501 Bielefeld}
\email{brussee@@mathematik.uni-bielefeld.de}

 \makeatletter

 \def\[{\begin{equation} \label}
 \def\]{\end{equation}}

 \theoremstyle{plain}

 \begingroup
 \theorembodyfont{\sl}

 \newtheorem{Theorem}{Theorem}
 \newtheorem{Corollary}[Theorem]{Corollary}

 \endgroup

 \theoremstyle{definition}

 \theoremstyle{remark}

 \newtheorem{Remark}[Theorem]{Remark}

 \newtheorem{Acknowledgment}{Acknowledgment}

 \def\eqalign#1{\mathalign\hfil\relax\relax\hfil#1\endmathalign}
 \def\mathalign#1#2#3#4{\null\,\vcenter\bgroup\let\\=\cr\openup\jot\m@th
 \ialign\bgroup\strut#1$\displaystyle##$#2&&#3$\displaystyle{}##$#4\crcr}
 \def\endmathalign{\crcr\egroup\egroup\,}

 \def\txt#1{\hbox{ #1 }}

\let\slant=/

\mathcode`\:="603A 
\mathchardef\:="303A 

\def\Q{{\Bbb Q}}
\def\Z{{\Bbb Z}}

\def\P{{\Bbb P}}

\def\numfrac#1#2{\mathchoice{{\textstyle{ #1\over#2}}}%
{{ #1\over#2}}{{#1/#2}}{{#1/#2}}}

\def\half{{\numfrac12}}

\let\next=\~
\let\~=\tilde
\let\tilde=\next

\let\next=\^
\let\^=\hat
\let\hat=\next

\def\"#1{{\accent"7F \if#1i\i\else#1\fi}}
\def\({\left(}
\def\){\right)}

\def\<#1>{\left<#1\right>}

\newif\ifcomment
\def\comment{\ifcomment}
\def\endcomment{}
\makeatother

\def\B#1{{}_{#1}B}
\def\S#1{{}_{#1}S}
\def\Pbar{{\bar{\P}}}
\def\ddt{{d \over dt}}
\def\D{{\mathcal D}}

\begin{document}

\maketitle

In this note we give a universal formula for the evaluation of the
Donaldson polynomials on $(-2)$-spheres, i.e. smooth spheres of
selfintersection $-2$.  Note that the blow-up formulas can be considered
as formulas for the evaluation of Donaldson polynomials on $(-1)$-spheres.
Such formulas come almost for free using Fintushel and Sterns method
to determine the ordinary blow-up formulas \cite{FS}. Fintushel and
Stern informed me that they found similar formulas with similar methods.
In fact comparing their formulas with my own ones, I realised that
remark \ref{failure} could be used to get clean formulas in the main theorem.

\begin{Acknowledgment}
Special  thanks to
both Fintushel and Stern for informing me about their results and
encouraging me to publish my own ones.
Further thanks to Stefan Bauer for getting me in the Bielefeld
vector bundle work group BVP, and Viktor Pidstrigatch and again  Stefan Bauer
for useful discussions.
\end{Acknowledgment}

Our conventions follow \cite{FS}.
We consider the (stable) Donaldson polynomial $D_c$ as a linear
form on $A(X) = S^\bullet H_*(X,\Q)$ determined by a class $c \in
H^2(X,\Z/4\Z)$. The generator of $H_0(X)$ is denoted $x$, and $\Q[x]$ is
considered as a subalgebra of $A(X)$.
We also use the convention that $D_c(a)$
means either evaluation of $D_c$ on $a \in
A(X)$ or the linear form defined by $D_c(a)(b) = D_c(ab)$
depending on the context.
Finally we extend $D_c$ linearly to the formal series
$A(X)[[t]]$ or the formal Laurent series $A(X)[[t]][t^{-1}]$
whenever convenient. Now the main result of \cite{FS} is that
there are  even respectively odd power series
$B(t)$, $S(t) \in \Q[x][[t]]$ such that for every $(-1)$-sphere
$e$ in a compact closed simply connected $4$-manifold with $b_+ >1$
\begin{eqnarray}
\label{*}
          D_c(e^{te}) &=   D_c(B(t)),
                \quad B(t)= \sigma_3(t) e^{-x\,t^2/6}
                &\qquad \txt{if}  c\cdot e \txt{is even}
\\
\label{**}
       D_c(e^{te}) &= D_{c-e}(S(t)),
              \quad S(t) = \sigma(t) e^{-x\,t^2/6}
              &\qquad \txt{if}  c\cdot e \txt{is odd}
\end{eqnarray}
on the subalgebra $A(e^\perp)$ generated by the orthogonal complement of
$e$ (this formulation is equivalent to the one where $c\cdot e = 0$).
Here $\sigma_3$ and $\sigma$ are certain quasi elliptic functions
with a power series expansion in $\Q[x][[t]]$.
Actually Fintushel and Stern mention that the simply connectedness
hypothesis is not necessary,  and the (simply connected) in our main
theorem below is meant to mean ``whatever they need''.

\begin{Theorem} \label{main}
Let $X$ be a (simply connected) closed compact oriented four manifold
with odd $b_+ >1$ containing a  $(-2)$ sphere $\tau$, then on
$A(\tau^\perp)$
\begin{eqnarray}
\label{maina}
       D_c(e^{t\tau}) &=
                      D_c(B^2(t))  +  D_{c +\tau}(S^2(t))
                                     & \txt{if} c\cdot\tau \txt{is even},
\\
\label{mainb}
         D_c(e^{t\tau})
                      &= D_c((BS' - B'S)(t)  +  \tau BS(t)) &\txt{if}
                                     c\cdot\tau \txt{is odd}.
\end{eqnarray}
\end{Theorem}

\comment
\begin{enumerate}
\item  There  are universal
even power series $\B2_0(t)$, $\B2_\tau$ with coefficients in  $\Q[x]$
such that if $c\cdot \tau \equiv 0 \pmod 2$
$$
       D_c(e^{t\tau}) = D_c(\cosh(t\tau))
                     = D_c(\B2_0(t)) + D_{c+\tau}(\B2_\tau(t))
       \quad\txt{on} (t)A(\tau^\perp)
$$
where $\B2_0$ and  $\B2_\tau$  are given in terms of $B$ and $S$ by
\begin{eqnarray}
\label{1a}
       \B2_0(t) &=&
       \sqrt{B(2t)}\, \cosh \(\half \int_0^{2t}{S\over B}(s)\, ds\),
\\
\label{1b}
       \B2_\tau(t) &=&
          \sqrt{B(2t)}\, \sinh\(\half \int_0^{2t}{S\over B}(s)\, ds \).
\end{eqnarray}
\item  There  are universal
even respectively odd power series $\S2_0(t)$, $\S2_1$
with coefficients in $\Q[x]$ such that if $c \cdot \tau \equiv 1 \pmod 2$
$$
       D_c(e^{t\tau}) = D_c(\S2_0(t) + \S2_1(t)\tau)
              \txt{on} A(\tau^\perp)
$$
where $\S2_0$ and  $\S2_1$ are given explicitly in terms of $B$ and $S$
by
\begin{eqnarray}
\label{2a}
         \S2_0(t) &=& e^{\half\int_0^{2t} {-B + S'\over S}(s)\,ds}
\\
\label{2b}
  \S2_1(t) &=&
   t \,e^{\half\int_0^{2t} {B + S' \over S}(s) - (2/ s) \,ds}
\end{eqnarray}
\end{enumerate}
\end{Theorem}

In the formulas above, divisions should be interpreted as follows.
The power series $B(t) = 1 + O(t^4)$,  hence $B$ is an invertible element
in $A(X)[[t]]$. Likewise $S(t) = t + O(t^3)$ is invertible
in the Laurent series $A(X)[[t]][t^{-1}]$.
\fi\endcomment

\begin{Remark}
Note that respectively Rubermann's relation \cite[Th 2.2]{FS}
\begin{equation}
\label{ruber}
         D_c(\tau^2) = 2D_{c + \tau} \qquad \txt{on} A(\tau^\perp), \quad
c\cdot \tau \txt{even}
\end{equation}
follows from  the formulas in the main theorem. In particular,
formula  \ref{maina} of the main  theorem implies
\begin{equation}
\label{main'}
         D_c(e^{t \tau}) = D_c(B^2(t) + S^2(t) \tau/2)
       \quad \txt{on} A(\tau^\perp), \quad c\cdot \tau \txt{even}.
\end{equation}
In this form it also implies Wieczorek's  relation \cite[Cor. 2.5]{FS}
\begin{equation}
\label{wieczorek}
         D_c(\tau^4) = D_c(-4 - 4x \tau^2).
\end{equation}
In fact to get these two relations
we only need the existence and form of the universal formulas for  $S$
and $B$ up to order 4. These two relations in turn determine the full
series for $S$ and $B$. Thus we will see  that the formulas for $B$ and
$S$ in \cite{FS} can be
proved from the magic formula  \eqref{***} for $(-3)$-spheres and the
expansion of $B$ and $S$  up to order 4 as gauge theoretic input (the
sign in \eqref{***} can easily be determined from $B$ and $S$ up to
order 4 instead of formula \eqref{ruber} as in \cite{FS}). Also
note that the evenness of $D_c(e^{t\tau})$ if $c \cdot \tau$ is even
is clear {\em a priori}
from the invariance of the Donaldson polynomial under the reflection in
$\tau$.
\end{Remark}

\comment
\begin{Remark}
Kronheimer and Mrowka obtain universal relations for tight
surfaces (i.e. those for which $2g(\Sigma) -2 = \Sigma^2$ c.f.
\cite[Def. 7.5]{KM}) under certain extra assumptions
on the selfintersection and genus. A $(-2)$-sphere is certainly tight,
but the extra assumptions are not satisfied. They ask
if their relations continue to hold without these
extra conditions \cite[section 9.2, question 2]{KM}.
For $(-2)$-spheres this is  not true as can be seen by comparing
formula \eqref{wieczorek} and  \cite[below equation 8.2]{KM}.
\end{Remark}
\fi\endcomment

The formulas simplify when $X$ is of simple type, i.e. if $D(x^2) = 4D$.
The elliptic functions $\sigma$ and $\sigma_3$ then degenerate  to
 respectively $\sinh$ and $\cosh$. The formulas are particularly
nice if we express them in terms of the Donaldson series
$$
         \D_c(h) = D_c((1+x/2)e^{h}), \quad h \in H_2(X)
$$
(this expression is a priori formal,  but in fact $\D_c$ is an entire real
analytic function  at least under the
assumption $b_1(X) = 0$, \cite{KM},\cite{FS:shape}).

\begin{Corollary} \label{simple type}
If in addition to the assumptions of the theorem, $X$ is of simple type
and $h \in \tau^\perp$, then
$$
\eqalign{
         \D_c(h + t\tau) &= e^{-t^2}(\cosh^2(t)\D_c(h) +
                 \sinh^2(t)\D_{c+\tau}(h)), & \txt{if}
                 c\cdot \tau \txt{is even},
\\
         \D_c(h + t\tau) &= e^{-t^2}(\D_c(h) + \sinh(2t)/2
\ddt|_{t=0}\D_c(h+t\tau)) & \txt{if} c \cdot \tau
                 \txt{is odd}.
}
$$
\end{Corollary}



\begin{pf}[theorem \ref{main}]
Theorem \ref{main} is an easy consequence of the  following magic
formula for  $(-3)$-spheres \cite[Theorem 2.4]{FS}. If $c\cdot \tau
\equiv 0 \pmod 2$ then
\begin{equation}
\label{***}      D_{c + \tau_3} = - D_c(\tau_3)   \txt{on} \tau_3^\perp
\end{equation}
As in \cite{FS} we blow-up to get relations for the $(-2)$-sphere.
In $X\#\Pbar^2$ the exceptional line is a $(-1)$-sphere $e$,
$\tau + e$ is represented by a $(-3)$-sphere,
$(\tau - 2e) \in (\tau + e)^\perp$, and $\tau \in e^\perp$.

First consider the case $c \cdot \tau$ even.
Then formula \eqref{***} gives the  identity
\begin{eqnarray*}
         D_{c + \tau + e}(e^{t(\tau - 2e)})
                 & \buildrel\eqref{**}\over=& D_{c+ \tau}(S(-2t) e^{t\tau})
\\
         & \buildrel\eqref{***}\over=& - D_c((\tau + e)e^{t(\tau -2e)})
\\
                 & \buildrel\eqref{*}\over=&
                 - D_c(B(-2t)\tau e^{t\tau} + B'(-2t)e^{t\tau}).
\end{eqnarray*}
on $\tau^\perp \cap e^\perp = \tau^\perp \subset H_2(X)$.
However, $B(t)$ is invertible in $\Q[x][[t]] \subset A(\tau^\perp)$.
Thus without loss of generality we can multiply by $1/B$ ``under $D_c$``.
With this remark and using the parity of $B$ and $S$ we get
\begin{equation}
\label{10}
         \ddt D_c(e^{t\tau}) = D_c({B'\over B}(2t) e^{t\tau})
                  + D_{c+\tau}({S\over B}(2t) e^{t\tau}).
\end{equation}
Now in this formula we can replace $c$ by $c+\tau$, and use that
$D_{c + 2\tau} = D_c$.
Hence if we define $D_{c,\pm} = D_c \pm D_{c+\tau}$, formula
\eqref{10} is equivalent to the two equations
\begin{equation}
\label{PMeqn}
         \ddt D_{c,\pm}(e^{t\tau}) =
                 D_{c,\pm}({B' \pm S \over B}(2t) e^{t\tau})
\end{equation}
This  differential equation determines $D_{c,\pm}(e^{t\tau})$
uniquely in terms of $D_{c,\pm}|_{\tau^\perp}$, i.e. there exists a
universal power series $\B2_\pm(t)$ with coefficients in $\Q[x]$ such that
$$
         D_{c,\pm}(e^{t\tau}) = D_{c,\pm}(\B2_\pm(t)) \txt{on}
            A(\tau^\perp).
$$
The series is only determined  modulo the intersection of the
kernels of $D_{c,\pm}(z)$ for all $z\in A(\tau^\perp))$,
but clearly  the following series is  a universal representative
\begin{equation}
\label{20}
         \B2_\pm(t) =
         e^{\int_0^t {B' \pm S \over B}(2s)\, ds} =
         \sqrt{B(2t)} e^{\pm{1 \over 2}\int_0^{2t} {S\over B}(s)\,ds}
        \in \Q[x][[t]].
\end{equation}
Formula \eqref{20} can then be rewritten to a universal formula
\begin{equation}
\label{univform}
         D_c(e^{t\tau}) = D_c(\B2_0(t)) + D_{c+\tau}(\B2_\tau(t))
                \txt{on} A(\tau^\perp).
\end{equation}
where $\B2_0 = \half(\B+ + \B-)$ and $\B2_\tau = \half(\B+ - \B-)$.

To get the clean formulas \eqref{maina} and \eqref{mainb}
in the theorem we simply apply the
universal formula \eqref{univform} to $\tau = e_1 + e_2 $ in $X \# 2
\Pbar^2$, and restrict to $H_2(X)$. For example, a small computation with
the blow-up formulas gives
$$
         D_c(S^2(t)) = D_{c + e_1 + e_2}(e^{t (e_1 + e_2)})|_{A(X)}
                    =  D_c(\B2_\tau)
$$
valid on all of $A(X)$. This proves the even case.

In case $c \cdot \tau$ is odd we proceed similarly.
We use the magic formula \eqref{***} for $D_{c+e}$ and $\tau_3 =
\tau + e$.
Now the reflection $R_\tau$ in $\tau$  is represented by an
orientation and
homology orientation preserving diffeomorphism,
so by the invariance of the polynomials.
$$
         R_{\tau}^*D_c = D_{c \pm \tau} = (-1)^{\tau^2} D_{c+\tau} =
                         D_{c+\tau}
$$
This granted,
\begin{eqnarray*}
         D_{c + \tau +2e}(e^{t(\tau - 2e)})
            &=& -D_{c+\tau}(B(-2t)(e^{t\tau})
\\
            &=& -D_c(B(2t)(e^{-t\tau})
\\
            &=& -D_{c+e}((\tau +e) e^{t(\tau -2e)}
\\
            &=& -D_c((-S(2t)\tau +S'(2t))e^{t\tau}).
\end{eqnarray*}
Some reworking then gives
\begin{eqnarray*}
         \ddt D_c(\cosh(t\tau)) &=& D_c({-B +S' \over S}(2t)  \cosh(t\tau))
\\
         \ddt D_c(\sinh(t\tau)) &=& D_c({ B+ S' \over S}(2t) \sinh(t\tau)).
\end{eqnarray*}

Since $B(0) = S'(0) = 1$, the first formula can be handled as above, and
we get a universal formula $D_c(\cosh(t\tau)) = D_c(\S2_0(t)$ on
$A(\tau^\perp)$ with an explicit expression
\begin{equation}\label{2S0}
         \S2_0 = e^{{1 \over 2}\int_0^{2t} {-B + S'\over S}(s)\,ds}.
\end{equation}
In the latter formula we have a differential equation with regular
singularities. It has a unique solution with
initial conditions  $D_c(\sinh(t\tau))|_{t = 0} = 0$,
$\ddt|_{t=0} D_c(\sinh(t\tau)) = D_c(\tau)$ which yields a universal
formula $D_c(\sinh(t\tau)) = D_c(\tau \S2_1(t))$ and an explicit
expression
$$
   \S2_1(t) = t \,e^{{1 \over 2}\int_0^{2t} {B + S' \over S}(s) - (2/ s) \,ds}
$$
Finally to get the clean formulas in the main theorem we blow-up twice
and compute $D_{c + e_1}(e^{t(e_1 + e_2)})|_{A(X)}$ and
$D_c(e^{t(e_1 + e_2)}(e_1 - e_2))|_{A(X)}$ in two different
ways as in the even
case.
\end{pf}

\begin{Remark} \label{failure}
It is by no means clear that say $\B2_0(t) = B^2(t)$ in $\Q[x][[t]]$,
i.e. without dividing out $\ker D_c$.  In fact, I hoped to get universal
identities for $D_c(x^n)$ like the simple type condition $D_c(x^2) = 4
D_c$ this way. Alas, expanding both sides by computer, I found the
series  agreed at least up to order 28 in $t$ for all series considered
here. It should be expected that $B^2(t) \pm S^2(t)$
satisfies the differential equation \eqref{PMeqn} before evaluation with
$D_c$

In the same spirit,
formulas for $(-2)$-spheres give us relations for $B$ and $S$ by
writing two orthogonal $(-1)$-spheres as a $(-2)$-sphere. For
example when  $c \in H^2(X)
\subset H^2(X \# 2 \Pbar^2)$  we get the identity
\begin{eqnarray}
\label{BB}
         D_c(B(u+v)B(u-v)) &=& D_c(B^2(u)B^2(v) - S^2(u) S^2(v))
\end{eqnarray}
on $A(X)$.
Likewise, writing a
$(-3)$-sphere as three $(-1)$-spheres and plugging this into formula
\eqref{***} gives identities like
\begin{eqnarray}
\label{BBB}
         D_c(S(u)S(v)S(u+v)) &=&
\\
\notag
         && \llap{$D_c(B'(u)B(v)B(u+v)$} + B(u)B'(v)B(u+v) - B(u)B(v)B'(u+v))
\end{eqnarray}
Fintushel and Stern informed me, that they reduced these identities to
identities between elliptic functions.
\end{Remark}

\section*{Appendix 1. Explicit series}

For the convenience of the reader we give the first few terms
of the  expansion of the relevant series:

\def\mon(#1){{t^{#1} \over #1 !} + }
$$
\eqalign{
B(t)&=
\\
&1  
+ (-2)\,\mon(4)
(8\,x)\mon(6)
\left (-4-32\,{x}^{2}\right )\mon(8)
\left (96\,x+128\,{x}^{3}\right )\mon(10)
\cr
&
\left (-408-960\,{x}^{2}-512\,{x}^{4}\right )\mon(12)
\left (7584\,x+7168\,{x}^{3}+2048\,{x}^{5}\right )\mon(14)
\cr
&
\left (13584-88320\,{x}^{2}-46080\,{x}^{4}-8192\,{x}^{6}\right )\mon(16)
%
%
%
\cdots
%
\cr
S(t) &=
\\
& t +  
(-x)\mon(3)
\left (2+{x}^{2}\right )\mon(5)
\left (-6\,x-{x}^{3}\right )\mon(7)
\left (-36+12\,{x}^{2}+{x}^{4}\right )\mon(9)
\cr
&
\left (564\,x-20\,{x}^{3}-{x}^{5}\right )\mon(11)
\left (-552-5916\,{x}^{2}+30\,{x}^{4}+{x}^{6}\right )\mon(13)
\cr
&
\left (+13848\,x+55404\,{x}^{3}-42\,{x}^{5}-{x}^{7}\right )
\mon(15)
%
%
\cdots
%
\cr
\B2_0(t) &= B^2(t)=
\cr
&
1 + 
(-4)\,\mon(4)
(16\,x)\mon(6)
\left (272-64\,{x}^{2}\right )\mon(8)
\left (-6528\,x+256\,{x}^{3}\right )\mon(10)
\cr
&
\left (7104 + 120576\,{x}^{2}-1024\,{x}^{4} \right )\mon(12)
\left (-561408\,x-2035712\,{x}^{3}+4096\,{x}^{5}\right )\mon(14)
%
%
%
%
%
%
\cdots
%
\cr
\B2_\tau(t) &= S^2(t) =
\cr
&2\,\mon(2)
(-8\,x)\mon(4)
\left (24+32\,{x}^{2}\right )\mon(6)
\left (-320\,x-128\,{x}^{3}\right )\mon(8)
\cr
&
\left (288+2688\,{x}^{2}+512\,{x}^{4}\right )\mon(10)
\left (10368\,x-18432\,{x}^{3}-2048\,{x}^{5}\right )\mon(12)
%
%
%
%
%
%
%
\cdots
%
%
\cr
\S2_0(t) &= (BS'- SB')(t) =
\\
&1 + 
(-x)\,\mon(2)
\left (8+{x}^{2}\right )\mon(4)
\left (-56\,x-{x}^{3}\right )\mon(6)
\left (-64+432\,{x}^{2}+{x}^{4}\right )\mon(8)
\cr
&
\left (-192\,x-3728\,{x}^{3}-{x}^{5}\right )\mon(10)
\left (26112+26688\,{x}^{2}+33272\,{x}^{4}+{x}^{6}\right )\mon(12)
%
%
%
%
%
%
\cdots
%
\cr
\S2_1(t)&= BS(t) =
\\
&
t +  
(-x)\mon(3)
\left (-8+{x}^{2}\right )\mon(5)
\left (120\,x-{x}^{3}\right )\mon(7)
\left (-576 -1200\,{x}^{2}+{x}^{4}\right )\mon(9)
\cr
&
\left (13632\,x+11024\,{x}^{3}-{x}^{5}\right )\mon(11)
\left (35328 -232896\,{x}^{2}-99576\,{x}^{4}+{x}^{6}\right )\mon(13)
%
%
%
%
%
%
\cdots
}
$$

\pagebreak

\end{document}